\documentclass{amsart}
\usepackage{times}
\usepackage{graphicx}

\numberwithin{equation}{section}

\usepackage{amssymb}

\DeclareMathOperator{\diam}{diam}

\newcommand{\sub}{\subseteq}
\newcommand{\sfrac}[2]{\mbox{$\frac{#1}{#2}$}}
\newcommand{\abs}[1]{\left| #1 \right|}
\newcommand{\norm}[1]{\left\| #1 \right\|}
\newcommand{\expr}[1]{\left( #1 \right)}
\newcommand{\set}[1]{\left\{ #1 \right\}}
\newcommand{\scalar}[1]{\left< #1 \right>}

\newcommand{\cD}{\mathcal{D}}
\newcommand{\cE}{\mathcal{E}}
\newcommand{\cS}{\mathcal{S}}
\newcommand{\cV}{\mathcal{V}}
\newcommand{\ex}{\mathbf{E}}
\newcommand{\pr}{\mathbf{P}}
\newcommand{\R}{\mathbf{R}}
\newcommand{\Z}{\mathbf{Z}}

\newtheorem{theorem}{Theorem}
\newtheorem{lemma}{Lemma}
\newtheorem{proposition}{Proposition}
\theoremstyle{definition}
\newtheorem{definition}{Definition}

\begin{document}

\title[Boundary Harnack inequality on the Sierpi{\'n}ski triangle]{Boundary Harnack inequality for $\alpha$-harmonic functions on the Sierpi{\'n}ski triangle}
\author{Kamil Kaleta}
\email{Kamil.Kaleta@pwr.wroc.pl}
\address{Institute of Mathematics and Computer Science \\ Wroc{\l}aw University of Technology \\ ul. Wybrze{\.z}e Wyspia{\'n}skiego 27 \\ 50-370 Wroc{\l}aw, Poland}
\author{Mateusz Kwa{\'s}nicki}
\email{Mateusz.Kwasnicki@pwr.wroc.pl}
\address{Institute of Mathematics and Computer Science \\ Wroc{\l}aw University of Technology \\ ul. Wybrze{\.z}e Wyspia{\'n}skiego 27 \\ 50-370 Wroc{\l}aw, Poland}

\begin{abstract}
We prove a uniform boundary Harnack inequality for nonnegative functions harmonic with respect to $\alpha$-stable process on the Sierpi{\'n}ski triangle, where $\alpha \in (0, 1)$. Our result requires no regularity assumptions on the domain of harmonicity.
\end{abstract}

\maketitle

\section{Introduction and main theorem}
\label{sec:intro}

The analysis and probability theory on fractals underwent rapid development in last twenty years, see~\cite{bib:b98, bib:dsv99, bib:strichartz99, bib:s06} and the references therein. Diffusion processes were constructed for the Sierpi{\'n}ski triangle~\cite{bib:bp88, bib:g87, bib:k87} and more generally for some simple nested fractals~\cite{bib:k01} and Sierpi{\'n}ski carpets~\cite{bib:bb89, bib:k00, bib:ks05, bib:pp00}. In~\cite{bib:s00} St{\'o}s introduced a class of \emph{subordinate} processes on $d$-sets, called \emph{$\alpha$-stable processes} on $d$-sets by analogy to the classical setting (see also~\cite{bib:k02}). Their nice scaling properties are similar to those of diffusion processes on $d$-sets, but their paths are no longer continuous. For the formal definition, see the Preliminaries section; here we only remark that in order to make the notion of $\alpha$-stability consistent with the scaling properties mentioned above, we depart from the notation of~\cite{bib:s00}. Namely, the $\alpha$-stable process below refers to the $(\frac{2 \alpha}{d_w})$-stable process in the sense of~\cite{bib:s00}. In particular, subordination yields $\alpha \in (0, d_w)$ rather than $\alpha \in (0, 2)$ as in~\cite{bib:s00}.

The theory of $\alpha$-stable processes on $d$-sets was further developed in~\cite{bib:bss03, bib:ck03, bib:k02}. In particular, it is known that the Harnack inequality holds true for nonnegative functions harmonic with respect to the $\alpha$-stable process ($\alpha$-harmonic functions) on a $d$-set $F$ whenever there is a diffusion process on $F$ and $\alpha \in (0, 1) \cup (d, d_w)$~\cite[Theorem~7.1]{bib:bss03}. Also, it is proved in~\cite[Theorem~8.6]{bib:bss03} that for the Sierpi{\'n}ski triangle a version of the boundary Harnack inequality holds for $\alpha \in (0, 1) \cup (d, d_w)$ if domain of harmonicity is a union of fundamental cells. The main result of this article extends this result for $\alpha \in (0, 1)$ to arbitrary open sets.

\begin{theorem}
\label{th:bhi}
Let $0 < \alpha < 1$. Let $B$ be the union of two adjacent cells of the infinite Sierpi{\'n}ski triangle $F$ with common vertex $x_0$, and let $B'$ be the union of the two twice smaller adjacent cells with common vertex $x_0$ (see Fig.~\ref{fig:bhp}).

There is a constant $c = c(\alpha)$ with the following property. Suppose that $D$ is an arbitrary open set in $F$. If $f$ and $g$ are nonnegative functions regular $\alpha$-harmonic in $D$ and vanishing on $D^c \cap B$, then
\begin{align}
\label{eq:bhi}
  \frac{f(x)}{g(x)}
  & \le
  c \, \frac{f(y)}{g(y)}
  \, , &
  x, y \in D \cap B'
  \, .
\end{align}
\end{theorem}

Our aim is to study the estimates and structure of $\alpha$-harmonic functions on $d$-sets. The present article is the case study of the Sierpi{\'n}ski triangle. It seems that the generalization of Theorem~\ref{th:bhi} to the case of more general simple nested fractals requires only minor changes in the proof, with the exception of the algebraic Lemma~\ref{lem:inequality}. Sierpi{\'n}ski carpets are another interesting fractals, studied e.g. in~\cite{bib:bgk, bib:hksz00, bib:s06:boundary}. However, proving an analogue of Theorem~\ref{th:bhi} in this case is much more difficult due to little knowledge about analysis on these sets.

Our argument follows the ideas of~\cite{bib:bkk08}, where isotropic $\alpha$-stable L{\'e}vy processes in $\R^d$ were considered. To adapt the argument for the fractal sets, two issues need to be resolved. First, a sufficiently smooth cutoff function is needed. In the case of Sierpi{\'n}ski triangle and some more general simple nested finite fractals, it can be constructed using \emph{splines}~\cite{bib:su00}. Second, the satisfactory estimate on the distribution of the process after it first exits from a ball is crucial for the proof of Lemma~\ref{lem:upper}. Such estimate is proved in~\cite[Lemmas~6.5 and~7.5]{bib:bss03} for general $d$-sets for $\alpha \in (0, 1)$ (this in fact is the only reason for the restriction on $\alpha$ in Theorem~\ref{th:bhi}). The problem whether similar result hold also for $\alpha \in [1, d_w)$ remains open.

The paper is organized as follows. In Preliminaries we recall the notions of Sierpi{\'n}ski triangle, fractional diffusion and $\alpha$-stable process, and construct cutoff functions on the Sierpi{\'n}ski triangle. Section~\ref{sec:harmonic} contains the proof of Theorem~\ref{th:bhi} and auxiliary lemmas.

\section{Preliminaries}
\label{sec:pre}

In this section we recall the construction of the unbounded Sierpi{\'n}ski triangle and the $\alpha$-stable process from~\cite{bib:bss03}, and collect some notation and facts.

\medskip

\noindent\textit{1. Sierpi{\'n}ski triangle.}
Let $F_0$ be the unit equilateral triangle, i.e. the closed triangle with vertices $p_1, p_2, p_3$, where $p_1 = (0, 0)$, $p_2 = (1, 0)$, $p_3 = (\frac{1}{2},\frac{\sqrt{3}}{2})$. Let $T_j$ denote the homothety with factor $\frac{1}{2}$ and center $p_j$, $j = 1, 2, 3$, and define a decreasing sequence of compact sets recursively by $F_{n+1} = T_1 F_n \cup T_2 F_n \cup T_3 F_n$. Let
\begin{align*}
  F_+
  & =
  \bigcap_{n=0}^\infty F_n
  \, .
\end{align*}
The set $F_+$ is the \emph{finite Sierpi{\'n}ski triangle}. Its mirror image about the vertical axis will be denoted by $F^-$. The \emph{infinite Sierpi{\'n}ski triangle} is defined by
\begin{align*}
  F
  & =
  \bigcup_{n=0}^\infty 2^n (F_+ \cup F_-)
  \, .
\end{align*}
For each $n \in \Z$, the infinite triangle $F$ is the union of the collection $\cS_n$ of uniquely determined isometric copies of $2^{-n} F_+$, called \emph{cells} of order $n$, or $n$-cells. The intersection of two distinct $n$-cells is either empty or contains a single point, called \emph{vertex} of order $n$, or $n$-vertex. In the latter case we say that the two $n$-cells are \emph{adjacent}. The set of $n$-vertices of $F$ is denoted by $\cV_n$. Two distinct $n$-vertices $u, v$ are \emph{adjacent}, $u \sim_n v$, if there is an $n$-cell containing both of them.

The infinite triangle $F$ is equipped with standard euclidean distance $\varrho(x, y) = |x - y|$; by $B(x, r)$ we denote an open ball in $F$. We remark that the intrinsic shortest-path metric $\varrho'$ is Lipschitz isomorphic to $\varrho$. Clearly, $F$ is arc-connected. Each $n$-cell contains three $n$-vertices, which constitute its topological boundary in $F$.

The Sierpi{\'n}ski triangle $F$ is a self-similar set of Hausdorff dimension $d = \frac{\log 3}{\log 2}$. If $D \sub F$ and $\diam D < \frac{1}{2}$, then $D$ is contained in some $0$-cell or in two adjacent $1$-cells. In either case there exist an open set $\tilde{D} \sub F_+ \cup F_-$ isometric to $D$. Furthermore, $F$ is invariant under homotheties with center at the origin and scale factor $2^n$, $n \in \Z$. These symmetries imply nice \emph{scaling properties} of various functions and measures on $F$.

Let $\mu$ denote the $d$-dimensional Hausdorff measure on $F$ so normalized that $\mu(F_+) = 1$. Any two isometric subsets of $F$ have equal measure, and $\mu(2^n E) = 3^n \mu(E)$ for any Borel $E \sub F$. For a function $f$ integrable on $E$,
\begin{align*}
  \int_E f(x) \mu(dx)
  & =
  3^{-n} \int_{2^n E} f(2^{-n} x) \mu(dx)
  \, .
\end{align*}

\smallskip

\noindent\textit{2. Calculus.}
In the past two decades calculus was developed for the \emph{finite} Sierpi{\'n}ski triangle, see e.g.~\cite{bib:b98, bib:fs92, bib:k01}. The extension to the infinite triangle is straightforward and we shall omit the details of this extension. Below we briefly introduce the concepts of the Laplace operator and the normal derivative.

\smallskip

Let $f$ be a continuous function on a cell $S \in \cS_n$. We define
\begin{align*}
  \cE_S(f, f)
  & =
  \lim_{k \rightarrow \infty} \expr{\sfrac{5}{3}}^k \sum_{u \sim_k w} \bigl(f(u) - f(w)\bigr)^2
  \, ,
\end{align*}
where the sum is taken over all pairs of neighbors $\set{u, w} \sub \cV_k \cap S$. We remark that the above limit is nondecreasing. Furthermore, for a continuous $f$ on $F$, we let
\begin{align*}
  \cE(f, f)
  & =
  \lim_{k \rightarrow \infty} \expr{\sfrac{5}{3}}^k \sum_{u \sim_k w} \bigl(f(u) - f(w)\bigr)^2
\end{align*}
with the summation over all neighbor pairs $\set{u, w} \sub \cV_k$. The domains of $\cE_S$ and $\cE$, denoted $\cD(\cE_s)$ and $\cD(\cE)$ respectively, consist of all functions $f$ for which the corresponding limits exist. Clearly $\cE(f, f) = \sum_{S \in \cS_n} \cE_S(f, f)$. Furthermore, $\cE_S$ and $\cE$ are regular local Dirichlet forms on $S$ and $F$, respectively~\cite{bib:fs92}. The Laplacian on $F$ is the self-adjoint (unbounded) operator on $L^2(F)$ associated to $\cE$; hence $\Delta f$ is the function in $L^2(F)$ satisfying
\begin{align*}
  \scalar{\Delta f, g}
  & =
  -\cE(f, g)
\end{align*}
for all $g \in \cD(\cE)$. The set of those $f$ for which $\Delta f$ exists is the domain of $\Delta$, denoted $\cD(\Delta)$. The Laplacian on a cell $S \in \cS_n$ is the operator $\Delta_S$ on $L^2(S)$ satisfying
\begin{align*}
  \scalar{\Delta_S f, g}
  & =
  - \cE_S(f, g)
\end{align*}
for all $g \in \cD(\cE_S)$ \emph{vanishing on the boundary of $S$}, with the domain $\cD(\Delta_S)$ being the set of $f \in \cD(\cE_S)$ for which such function exists. We emphasize that $\Delta_S$ is not a self-adjoint operator, and $\cD(\Delta_S)$ is larger than the domains of Dirichlet or Neumann Laplacians on $S$, see~\cite{bib:s06}.

\smallskip

Let $v \in \cV_n$ be a vertex of an $n$-cell $S$. For each $k \ge n$ there is a unique $k$-cell $S_k \sub S$ such that $v \in S_k$. Let $u_k, w_k$ denote the other $k$-vertices of $S_k$. The \emph{(outer) normal derivative} for $S$ and a function $f : S \rightarrow \R$ is defined by
\begin{align*}
  \partial_S f(v)
  & =
  \lim_{k \rightarrow \infty} \expr{\sfrac{5}{3}}^k \bigl(2 f(v) - f(u_k) - f(w_k)\bigr) \, ,
\end{align*}
provided the limit exists. Clearly, at each $n$-vertex $v$ there exist two normal derivatives $\partial_{S_1} f(v)$ and $\partial_{S_2} f(v)$, for the two adjacent $n$-cells $S_1$ and $S_2$ with common vertex $v$. For $f$ in the domain of $\Delta$ (or $\Delta_S$ with any $S \supseteq S_1 \cup S_2$) both $\partial_{S_1} f(v)$ and $\partial_{S_2} f(v)$ exist~\cite{bib:k01, bib:s06} and $\partial_{S_1} f(v) + \partial_{S_2} f(v) = 0$~\cite{bib:su00}. Furthermore, for $f \in \cD(\Delta_S)$ and $g \in \cD(\cE_S)$ we have by~\cite{bib:k01}
\begin{align*}
  \cE_S(f, g)
  & =
  - \scalar{\Delta_S f, g} + \sum_{v \in \partial S} \partial_S f(v) \, g(v)
  \, .
\end{align*}
For a more detailed introduction to the topic, the reader is referred to e.g.~\cite{bib:b98, bib:fs92, bib:k01, bib:s06, bib:t98}.

\medskip

\noindent\textit{3. Diffusion and stable processes.}
There exists a fractional diffusion on $F$ \cite{bib:b98, bib:bp88}. That is, there is a Feller diffusion $(Z_t)$ with state space $F$, such that its transition density function $q_t(x, y)$ (with respect to the Hausdorff measure $\mu$) is jointly continuous in $(x, y) \in F \times F$ for every $t > 0$ and satisfies
\begin{align}
\label{eq:diffusion}
  \frac{c_1'}{t^{\frac{d}{d_w}}} \exp \expr{-c_2' \frac{\varrho(x, y)^{\frac{d_w}{d_w - 1}}}{t^{\frac{1}{d_w - 1}}}}
  \le
  q_t(x, y)
  & \le
  \frac{c_1}{t^{\frac{d}{d_w}}} \exp \expr{-c_2 \frac{\varrho(x, y)^{\frac{d_w}{d_w - 1}}}{t^{\frac{1}{d_w - 1}}}}
\end{align}
for some positive $c_1, c_1', c_2, c_2'$ and all $t > 0$, $x, y \in F$. The constant $d_w = \frac{\log 5}{\log 2} \approx 2.322$ is the \emph{walk dimension} of $F$. This diffusion corresponds to the Dirichlet form $\cE$. We remark that the process corresponding to the Dirichlet form $\cE_S$ can be viewed as $(Z_t)$ reflected at the boundary of $S$.

By $Q_t$ we denote the transition operators of $(Z_t)$, $Q_t f(x) = \int q_t(x, y) f(y) \mu(dy)$, acting on either $L^2(F)$ or $C_0(F)$. The infinitesimal generator of $Q_t$ acting on $L^2(F)$ is precisely the Laplacian $\Delta$ defined in the previous paragraph. It agrees with the infinitesimal generator on $C_0(F)$ on the intersection of domains. The probability measure of the process $Z_t$ starting at $x \in F$ is denoted by $\pr^x$, and the corresponding expected value by $\ex^x$.

\smallskip

We fix $\alpha \in (0, d_w)$. Let $(Y_t)$ be the strictly $(\frac{\alpha}{d_w})$-stable subordinator, i.e. the nonnegative L{\'e}vy process on $\R$ with Laplace exponent $u^{\alpha/d_w}$~\cite{bib:b96, bib:bg68, bib:sato99}. We assume that $(Y_t)$ and $(Z_t)$ are stochastically independent. The subordinate process $(X_t)$, defined by $X_t = Z(Y_t)$, will be called the $\alpha$-stable process on $F$~\cite{bib:s00}. In~\cite{bib:s00}, $(X_t)$ is called $(\frac{2 \alpha}{d_w})$-stable; the change in notation is motivated by the scaling properties indicated below.

If $\eta_t(u)$ denotes the transition density of $(Y_t)$, then
\begin{align*}
  p_t(x, y)
  & =
  \int_0^\infty q_u(x, y) \, \eta_t(u) \, du
\end{align*}
defines the transition density of $(X_t)$. The corresponding transition operators $P_t f(x) = \int p_t(x, y) f(y) \mu(dy)$ form a semigroup on $C_0(F)$ and on $L^2(F)$, and the $L^2(F)$ infinitesimal generator of this semigroup is
\begin{align*}
  -(-\Delta)^{\frac{\alpha}{d_w}} f(x)
  =
  \lim_{t \searrow 0} \frac{P_t f(x) - f(x)}{t}
  \, ;
\end{align*}
the fractional power here is understood in the sense of spectral theory of unbounded operators on $L^2(F)$. We remark that both $q_t(x, y)$ and $p_t(x, y)$ have nice scaling properties,
\begin{align*}
  q_t(x, y)
  & =
  3^n q(2^{d_w n} t, 2^n x, 2^n y)
  \, , &
  p_t(x, y)
  & =
  3^n p(2^{\alpha n} t, 2^n x, 2^n y)
  \, .
\end{align*}
Furthermore, the $\pr^x$ law of $(X_t)$ is equal to the $\pr^{2^n x}$ law of $(2^{-n} \, X_{2^{\alpha n} t})$, and a similar relation holds for $Z_t$ with $\alpha$ substituted by $d_w$. If $f_n(x) = f(2^{-n} x)$, then for suitable $f$ we also have $\Delta f(x) = 2^{d_w n} \Delta f_n(2^n x)$ and $(-\Delta)^{\alpha/d_w} f(x) = 2^{\alpha n} (-\Delta)^{\alpha/d_w} f_n(2^n x)$.

There is $c_3 > 0$ such that~\cite[Theorem 37.1]{bib:d74}
\begin{align*}
  \lim_{u \rightarrow \infty} u^{1 + \frac{\alpha}{d_w}} \eta_1(u)
  & =
  \frac{\alpha}{2 \Gamma(1-\frac{\alpha}{d_w})}
  \, ,
  &
  \eta_1(u)
  & \le
  c_3 \min(1, u^{-1 - \frac{\alpha}{d_w}})
  \, , && u > 0
  \, .
\end{align*}
Denote $A_\alpha = \alpha / (2 \Gamma(1 - \frac{\alpha}{d_w}))$. By the scaling property,
\begin{align*}
  \eta_t(u)
  & =
  t^{-\frac{d_w}{\alpha}} \eta_1(t^{-\frac{d_w}{\alpha}} u)
  \, ,
  & t, u > 0
  \, ,
\end{align*}
we have
\begin{align}
\label{eq:eta}
  \lim_{t \searrow 0} \frac{\eta_t(u)}{t}
  & =
  A_\alpha u^{-1-\frac{\alpha}{d_w}}
  \, ,
  &
  \frac{\eta_t(u)}{t}
  & \le
  c_3 \min(t^{-\frac{d_w}{\alpha}}, t u^{-1 - \frac{\alpha}{d_w}})
  \, , && u > 0
  \, .
\end{align}
This formula will be used in Lemma~\ref{lem:spline}. We remark that~\eqref{eq:eta} and~\eqref{eq:diffusion} yield estimates of $p_t(x, y)$, see~\cite{bib:bss03}.

\smallskip

For a (relatively) open $D \sub F$, the first exit time of $D$,
\begin{align*}
  \tau_D
  & =
  \inf \set{t \ge 0 \; : \; X_t \notin D}
  \, ,
\end{align*}
is the stopping time. If $D$ is bounded then $\tau_D < \infty$ a.s., and the Green operator,
\begin{align*}
  G_D f(x)
  & =
  \ex^x \int_0^{\tau_D} f(X_t) \, dt
  \, ,
\end{align*}
has a nonnegative symmetric kernel $G_D(x, y)$ jointly continuous in $(x, y) \in D \times D$, and integrable in $y \in D$ for all $x \in D$~\cite[Section~5]{bib:bss03}. In particular, $G_D$ is a bounded operator on $C(D)$ and on $L^\infty(D)$, and $G_D f(x) \le \norm{f}_\infty \ex^x \tau_D$.

\smallskip

If $B$ is an open set such that $\overline{B}$ is compact and $\overline{B} \sub D$, we write $B \Subset D$. 

\begin{definition}
\label{def:harmonic}
A function $f : F \rightarrow [0, \infty)$ is \emph{$\alpha$-harmonic} in open $D \sub F$ if
\begin{align}
\label{eq:harmonic}
  f(x)
  & =
  \ex^x f(X(\tau_B))
  &
  \text{for every open $B \Subset D$ and $x \in B$.}
\end{align}
If~\eqref{eq:harmonic} holds for all $B \sub D$ (in particular for $B = D$) then $f$ is \emph{regular $\alpha$-harmonic} in $D$.
\end{definition}

By the strong Markov property, if $f(x) = \ex^x g(X(\tau_D))$ for some nonnegative $g$, then $f$ is regular $\alpha$-harmonic in $D$.

If $f$ is (regular) $\alpha$-harmonic in $D$, then $f_n(x) = f(2^{-n} x)$ is (regular) $\alpha$-harmonic in $2^n D$. Furthermore,
\begin{align*}
  \ex^x \tau_D
  & =
  2^{-\alpha n} \, \ex^{2^n x} (\tau_{2^n D})
  \, .
\end{align*}
We will use these and similar scaling properties without explicit reference.

\medskip

\noindent\textit{4. Splines.}
To construct a sufficiently smooth cutoff function $\varphi$ we will use the concept of \emph{splines} on the Sierpi{\'n}ski triangle~\cite{bib:su00}. First we prove some simple properties of a certain function on a cell of $F$.

Fix $S \in \cS_n$ and let $v_1, v_2, v_3$ be its vertices. Let $\varphi_0$ denote the function $f^{(1)}_{01}$ of~\cite{bib:su00}, the element of the \emph{better basis}, rescaled to $S$. This is a biharmonic function on $S$ (i.e. $(\Delta_S)^2 \varphi_0 = 0$) satisfying $\varphi_0(v_1) = 1$, $\varphi_0(v_2) = \varphi_0(v_3) = 0$ and $\partial_S \varphi_0(v_j) = 0$ for $j = 1, 2, 3$.

\begin{proposition}
\label{prop:recurrence}
Suppose that $u_1, u_2, u_3$ are vertices of a $k$-cell $S' \sub S$. Let $w_1 = u_1$, and let $w_2, w_3$ be the midpoints of line segments $u_1 u_2$ and $u_1 u_3$, respectively. The three points $w_1, w_2, w_3$ are vertices of a $(n+1)$-cell $S'' \sub S'$. We have
\begin{align*}
  \begin{pmatrix}
    \varphi_0(w_1) \\
    \varphi_0(w_2) \\
    \varphi_0(w_3) \\
    \expr{\frac{3}{5}}^{k+1} \partial_{S''} \varphi_0(w_1) \\
    \expr{\frac{3}{5}}^{k+1} \partial_{S''} \varphi_0(w_2) \\
    \expr{\frac{3}{5}}^{k+1} \partial_{S''} \varphi_0(w_3)
  \end{pmatrix}
  & =
  \frac{1}{75}
  \begin{pmatrix}
    75 & 0 & 0 & 0 & 0 & 0 \\
    36 & 36 & 3 & -7 & -7 & -1 \\
    36 & 3 & 36 & -7 & -1 & -7 \\
    0 & 0 & 0 & 45 & 0 & 0 \\
    -90 & 90 & 0 & 15 & -15 & 0 \\
    -90 & 0 & 90 & 15 & 0 & -15
  \end{pmatrix}
  \begin{pmatrix}
    \varphi_0(u_1) \\
    \varphi_0(u_2) \\
    \varphi_0(u_3) \\
    \expr{\frac{3}{5}}^k \partial_{S'} \varphi_0(u_1) \\
    \expr{\frac{3}{5}}^k \partial_{S'} \varphi_0(u_2) \\
    \expr{\frac{3}{5}}^k \partial_{S'} \varphi_0(u_3)
  \end{pmatrix}
  \, .
\end{align*}
\end{proposition}

\begin{proof}
Formula~(5.8) of~\cite{bib:su00} states that
\begin{align*}
  \begin{pmatrix}
    \varphi_0(w_1) \\
    \varphi_0(w_2) \\
    \varphi_0(w_3) \\
    \expr{\frac{1}{5}}^{k+1} \Delta \varphi_0(w_1) \\
    \expr{\frac{1}{5}}^{k+1} \Delta \varphi_0(w_2) \\
    \expr{\frac{1}{5}}^{k+1} \Delta \varphi_0(w_3)
  \end{pmatrix}
  & =
  \frac{1}{25}
  \begin{pmatrix}
    25 & 0 & 0 & 0 & 0 & 0 \\
    10 & 10 & 5 & -3/5 & -3/5 & -7/9 \\
    10 & 5 & 10 & -3/5 & -7/9 & -3/5 \\
    0 & 0 & 0 & 5 & 0 & 0 \\
    0 & 0 & 0 & 2 & 2 & 1 \\
    0 & 0 & 0 & 2 & 1 & 2
  \end{pmatrix}
  \begin{pmatrix}
    \varphi_0(u_1) \\
    \varphi_0(u_2) \\
    \varphi_0(u_3) \\
    \expr{\frac{1}{5}}^k \Delta \varphi_0(u_1) \\
    \expr{\frac{1}{5}}^k \Delta \varphi_0(u_2) \\
    \expr{\frac{1}{5}}^k \Delta \varphi_0(u_3)
  \end{pmatrix}
  \, .
\end{align*}
For brevity, we write this formula as $\mathbf{d}_{k+1} = A \, \mathbf{d}_k$. Furthermore, by~(3.5) and~(5.9) of~\cite{bib:su00}, scaling and the construction of $\varphi_0$,
\begin{align*}
  \begin{pmatrix}
    \varphi_0(u_1) \\
    \varphi_0(u_2) \\
    \varphi_0(u_3) \\
    \expr{\frac{1}{5}}^k \Delta \varphi_0(u_1) \\
    \expr{\frac{1}{5}}^k \Delta \varphi_0(u_2) \\
    \expr{\frac{1}{5}}^k \Delta \varphi_0(u_3)
  \end{pmatrix}
  =
  \begin{pmatrix}
    1 & 0 & 0 & 0 & 0 & 0 \\
    0 & 1 & 0 & 0 & 0 & 0 \\
    0 & 0 & 1 & 0 & 0 & 0 \\
    -30 & 15 & 15 & 11 & -4 & -4 \\
    15 & -30 & 15 & -4 & 11 & -4 \\
    15 & 15 & -30 & -4 & -4 & 11
  \end{pmatrix}
  \begin{pmatrix}
    \varphi_0(u_1) \\
    \varphi_0(u_2) \\
    \varphi_0(u_3) \\
    \expr{\frac{3}{5}}^k \partial_{S'} \varphi_0(u_1) \\
    \expr{\frac{3}{5}}^k \partial_{S'} \varphi_0(u_2) \\
    \expr{\frac{3}{5}}^k \partial_{S'} \varphi_0(u_3)
  \end{pmatrix}
  \, .
\end{align*}
Again, this can be written in short as $\mathbf{d}_k = B \, \mathbf{c}_k$. A similar formula holds for $w_1, w_2, w_3$ and $S''$, in symbols: $\mathbf{d}_{k+1} = B \, \mathbf{c}_{k+1}$. It follows that $\mathbf{c}_{k+1} = B^{-1} A \, B \, \mathbf{c}_k$ and the proposition follows.
\end{proof}

\begin{lemma}
\label{lem:inequality}
The function $\varphi_0$ satisfies $0 \le \varphi_0(x) \le 1$ for all $x \in S$.
\end{lemma}

\begin{proof}
Let $S' \sub S$ be a $k$-cell with vertices $u_1, u_2, u_3$. We consider the following condition:
\begin{align}
\label{eq:splineok}
  \varphi_0(u_j) & \ge 0
  \quad \text{and} \quad
  \expr{\sfrac{3}{5}}^k \abs{\partial_{S'} \varphi_0(u_j)} \le 3 \varphi_0(u_j)
  \, , &
  j & \in \set{1, 2, 3}
  \, .
\end{align}
Note that~\eqref{eq:splineok} holds when $S' = S$. We claim that if~\eqref{eq:splineok} is satisfied for a $k$-cell $S' \sub S$, then it holds for each of the three $(k+1)$-cells $S'' \sub S'$.

Indeed, assume~\eqref{eq:splineok}, and let $S'', w_1, w_2, w_3$ be defined as above. By Proposition~\ref{prop:recurrence},
\begin{align*}
  \varphi_0(w_2)
  & =
  \expr{\sfrac{12}{25} \, \varphi_0(u_1) - \sfrac{7}{75} \expr{\sfrac{3}{5}}^k \partial_{S'} \varphi_0(u_1)}
  \\ & \quad +
  \expr{\sfrac{12}{25} \, \varphi_0(u_2) - \sfrac{7}{75} \expr{\sfrac{3}{5}}^k \partial_{S'} \varphi_0(u_2)}
  \\ & \quad +
  \expr{\sfrac{1}{25} \, \varphi_0(u_3) - \sfrac{1}{75} \expr{\sfrac{3}{5}}^k \partial_{S'} \varphi_0(u_3)}
  \ge
  0
  \, ,
\end{align*}
and
\begin{align*}
  3 \varphi_0(w_2) - \expr{\sfrac{3}{5}}^{k+1} \partial_{S''} \varphi_0(w_2)
  & =
  \expr{\sfrac{66}{25} \, \varphi_0(u_1) - \sfrac{12}{25} \expr{\sfrac{3}{5}}^k \partial_{S'} \varphi_0(u_1)}
  \\ & \quad +
  \expr{\sfrac{6}{25} \, \varphi_0(u_2) - \sfrac{2}{25} \expr{\sfrac{3}{5}}^k \partial_{S'} \varphi_0(u_2)}
  \\ & \quad +
  \expr{\sfrac{3}{25} \, \varphi_0(u_3) - \sfrac{1}{25} \expr{\sfrac{3}{5}}^k \partial_{S'} \varphi_0(u_3)}
  \ge 0
  \, .
\end{align*}
A similar calculation shows that $3 \varphi_0(w_2) + \expr{\frac{3}{5}}^{k+1} \partial_{S''} \varphi_0(w_2) \ge 0$. By symmetry, similar formulas hold also for $w_3$. This proves our claim.

By induction, $\varphi_0$ is nonnegative on every vertex in $S$. By continuity, $\varphi_0 \ge 0$ everywhere on $S$. Finally, the function $(1 - \varphi_0)$ is the sum of two copies of $\varphi_0$ with $v_1, v_2, v_3$ arranged in a different order, and this yields the inequality $\varphi_0 \le 1$.
\end{proof}

Suppose that a finite set of cells $\cS \sub \cS_n$ is given. Let $V \sub \cV_n$ be the set of vertices of cells from $\cS$. Define the cutoff function $\varphi$ in the following way. On each $n$-cell $S \in \cS_n$ with vertices $v_1, v_2, v_3$ we let:
\begin{align*}
  \varphi & = 1 & \text{if } & v_1, v_2, v_3 \in V \, ; \\
  \varphi & = 0 & \text{if } & v_1, v_2, v_3 \notin V \, ; \\
  \varphi & = \varphi_0 & \text{if } & v_1 \in V , \, v_2, v_3 \notin V \, ; \\
  \varphi & = 1 - \varphi_0 & \text{if } & v_1 \notin V , \, v_2, v_3 \in V \, .
\end{align*}
Here $v_1, v_2, v_3$ are arranged in a suitable order, so that one of the above conditions is satisfied. Observe $\varphi = 1$ on each $n$-cell in $\cS$ and $\varphi = 0$ on each $n$-cell disjoint with all cells from $\cS$. Furthermore, the definition of $\varphi$ is consistent in the following sense. When $v \in \cV_n$ is a common vertex of two $n$-cells $S_1, S_2$, then $\varphi$ is continuous at $v$ (that is, the definitions of $\varphi$ on $S_1$ and $S_2$ agree at $v$), and $\partial_{S_1} \varphi(v) + \partial_{S_2} \varphi(v) = 0$ because both normal derivatives vanish. Hence, by Theorem~4.4 of~\cite{bib:su00} (extended to the infinite triangle), $\varphi$ belongs to domain of $\Delta$, and $\Delta \varphi$ is essentially bounded on $F$. This smoothness property of $\varphi$ is used in the following result.

\begin{lemma}
\label{lem:spline}
The function $\varphi$ defined above belongs to the $C_0(F)$-domain of $-(-\Delta)^{\alpha/d_w}$.
\end{lemma}

\begin{proof}
By the Fubini Theorem,
\begin{align*}
  \frac{P_t \varphi(x) - \varphi(x)}{t}
  & =
  \frac{1}{t} \expr{\int \expr{\int_0^\infty \eta_t(u) \, q_u(x, y) \, du} \varphi(y) \, \mu(dy) - \varphi(x)}
  \\ & =
  \frac{1}{t} \int_0^\infty \eta_t(u) \expr{\int \varphi(y) q_u(x, y) \mu(dy) - \varphi(x)} du
  \\ & =
  \int_0^\infty \frac{\eta_{t}(u)}{t} \expr{Q_u \varphi(x) - \varphi(x)} du
  \, .
\end{align*}
We will show that
\begin{align}
\label{eq:convergence}
  \frac{P_t \varphi(x) - \varphi(x)}{t}
  & \rightarrow
  A_\alpha \int_0^\infty u^{-1 - \frac{\alpha}{d_w}} \expr{Q_u \varphi(x) - \varphi(x)} du
\end{align}
in the supremum norm. We have
\begin{align*}
  \abs{Q_u\varphi(x)-\varphi(x)}
  & =
  \abs{\int_0^u \frac{d}{ds} Q_s \varphi(x) ds}
  =
  \abs{\int_0^u Q_s \Delta \varphi(x) ds}
  \le
  u \norm{\Delta \varphi(x)}_\infty
  \, .
\end{align*}
Furthermore, $|Q_u \varphi(x) - \varphi(x)| \le 2 \, \|\varphi(x)\|_\infty$. It follows that
\begin{align*}
  &
  \abs{\frac{P_t \varphi(x) - \varphi(x)}{t} - A_\alpha \int_0^\infty u^{-1 - \frac{\alpha}{d_w}} (Q_u \varphi(x) - \varphi(x)) du}
  \\ & \quad \le
  \int_0^\infty \abs{\frac{\eta_t(u)}{t} - A_\alpha u^{-1 - \frac{\alpha}{d_w}}} \abs{Q_u \varphi(x) - \varphi(x)} du
  \\ & \quad \le
  \int_0^\infty \abs{\frac{\eta_t(u)}{t} - A_\alpha u^{-1 - \frac{\alpha}{d_w}}} \min \expr{u \norm{\Delta \varphi(x)}_\infty, 2 \norm{\varphi(x)}_\infty} du
  \, .
\end{align*}
By~\eqref{eq:eta} and dominated convergence theorem, the above integral converges to zero as $t \rightarrow 0$. This proves the uniform convergence in~\eqref{eq:convergence}, and the lemma follows.
\end{proof}

\section{Estimates of $\alpha$-harmonic functions}
\label{sec:harmonic}

We generally follow the proof of Theorem~1 of~\cite{bib:bkk08}. The argument incorporates some ideas from earlier works, particularly~\cite{bib:b97, bib:sw99}.

\begin{lemma}
\label{lem:escape}
For every $p_1, p_2$ such that $0 < p_1 < p_2 \le 1$ there is a constant $c_4 = c_4(p_1, p_2, \alpha)$ such that if $D$ is an open subset of the ball $B(v, p_2 2^{-m})$ for some $v \in \cV_m$, then
\begin{align}
\label{eq:escape}
  \pr^x(X(\tau_D) \notin B(v, p_2 2^{-m}))
  & \le
  c_4 \, 2^{\alpha m} \, \ex^x \tau_D
  \, , &
  x \in D \cap B(v, p_1 2^{-m})
  \, .
\end{align}
\end{lemma}

\begin{proof}
Note that formula~\eqref{eq:escape} is invariant under homothety with center $0$ and scale factor $2^m$. Hence we may assume that $m = 0$.

Choose $n$ large enough, so that any two $n$-cells $S, S'$ with $S \cap B(v, p_1) \ne \varnothing$ and $S' \cap B(v, p_2)^c \ne \varnothing$ have no common vertex. Let $V \sub \cV_n$ be the set of all vertices $v$ of $n$-cells $S$ satisfying $S \cap B(v, p_1) \ne \varnothing$. Let $\varphi$ be the cutoff function corresponding to $V$ constructed in the previous section. Clearly $\varphi = 1$ on $B(v, p_1)$ and $\varphi = 0$ on $(B(v, p_2))^c$. By Lemmas~\ref{lem:inequality} and~\ref{lem:spline}, $0 \le \varphi \le 1$ and $(-\Delta)^{\alpha/d_w} \varphi$ is essentially bounded. By formula~(5.8) of~\cite{bib:d65}, for $x \in D \cap B(v, p_1)$ we have
\begin{align*}
  \pr^x(X(\tau_D) \notin B(v, p_2))
  & =
  \ex^x \expr{\varphi(x) - \varphi(X(\tau_D)) \; ; \; X(\tau_D) \notin B(v, p_2)}
  \\ & \le
  \ex^x \expr{\varphi(x) - \varphi(X(\tau_D))}
  \\ & =
  \ex^x \int_0^{\tau_D} (-\Delta)^{\frac{\alpha}{d_w}} \varphi(X_t) \, dt
  \\ & \le
  \norm{(-\Delta)^{\alpha/d_w} \varphi}_\infty \ex^x \tau_D
  \, .
\qedhere
\end{align*}
\end{proof}

The proof of the next lemma hinges on the following two results of~\cite{bib:bss03}. For some positive $c_5 = c_5(\alpha)$ and $c_5' = c_5'(\alpha)$,
\begin{equation}
\label{eq:iw}
\begin{split}
  c_5' \int_{\overline{D}^c} \int_D \frac{G_D(x, y) \, f(z)}{\varrho(y, z)^{d + \alpha}} \, \mu(dy) \, \mu(dz)
  & \le
  \ex^x f(X(\tau_D))
  \\ & \le
  c_5 \int_{\overline{D}^c} \int_D \frac{G_D(x, y) \, f(z)}{\varrho(y, z)^{d + \alpha}} \, \mu(dy) \, \mu(dz)
\end{split}
\end{equation}
for all nonnegative $f$ with $f(z) = 0$ for $z \in \partial D$~\cite[Corollary 6.2]{bib:bss03}.

Suppose that $\alpha \in (0, 1)$. From the proof of Theorem~7.1 in the transient case in~\cite[Section~7.2]{bib:bss03} it follows that given any $v \in F$ and $r > s > 0$, there is a kernel function $P_{v, r, s}(x, y)$, $x \in B(v, s)$, $y \in (B(v, s))^c$ with the following two properties. There is $c_6 = c_6(\alpha, q)$ such that $P_{v, r, q}(x, y) \le c_6 r^{-d}$ for all $x, y$. Whenever $f$ is regular $\alpha$-harmonic in $B(v, r)$,
\begin{align*}
  f(x)
  & =
  \int P_{v, r, s}(x, y) \, f(y) \, \mu(dy)
  \, , &
  x \in B(v, s)
  \, .
\end{align*}
When $s = \frac{1}{4} \, r$, $P_{v, r, s}$ equals $2 r^{-1} P$, where $P$ is the function defined in~\cite{bib:bss03} with twice smaller $r$. For a general $s \in (0, r)$, $P_{v, r, s}$ is defined in a similar manner by changing the integration range $[r, 2r]$ to $[2s', 2r]$ in the definition of $P$, with some $s' \in (s, r)$, see~\cite{bib:bss03}.

If $D \sub B$, then by the strong Markov property,
\begin{align*}
  \pr^x(X(\tau_D) \in E)
  & \le
  \pr^x(X(\tau_B) \in E)
  \, , &
  \text{if $E \sub B^c$}
  \, .
\end{align*}
This \emph{monotonicity} of exit distributions implies that if $f$ is regular $\alpha$-harmonic in $D$ and $f = 0$ on $B \setminus D$, then
\begin{align*}
  f(x)
  & \le
  \ex^x f(X(\tau_B))
  \, , &
  x & \in D
  \, .
\end{align*}
This and the construction of $P_{v, r, s}$ yields that if $f$ is regular $\alpha$-harmonic in an open $D \sub B(v, r)$, then
\begin{align}
\label{eq:regularization}
  f(x)
  & \le
  \int P_{v, r, s}(x, y) \, f(y) \, \mu(dy)
  \, , &
  x \in B(v, s)
  \, .
\end{align}

\smallskip

For $v \in F$, $r > 0$ and a nonnegative function $f$ define
\begin{align*}
  \Lambda_{v, r}(f)
  & =
  \int_{B(v, r)^c} \varrho(y, v)^{-d - \alpha} f(y) \, \mu(dy)
  \, .
\end{align*}
Observe that $\Lambda_{v, r}(f) = 2^{\alpha n} \Lambda_{2^n v, 2^n r}(f_n)$, where $f_n(x) = f(2^{-n} x)$.

\begin{lemma}
\label{lem:upper}
Suppose that $\alpha \in (0, 1)$ and $0 < p_3 < p_5 \le 1$. There is a constant $c_7 = c_7(p_3, p_5, \alpha)$ with the following property. If a nonnegative function $f$ is regular $\alpha$-harmonic in an open $D \sub B(v, p_5 2^{-m})$, where $v \in \cV_m$, and vanishes on $D^c \cap B(v, p_5 2^{-m})$, then
\begin{align}
\label{eq:upper}
  f(x)
  & \leq
  c_7 \, 2^{-\alpha m} \Lambda_{v, p_3 2^{-m}}(f)
  \, , &
  x \in D \cap B(v, p_3 2^{-m})
  \, .
\end{align}
\end{lemma}

\begin{proof}
As in the previous lemma, \eqref{eq:upper} is invariant under dilations and hence we may assume that $m = 0$. Fix any $p_4$ such that $p_3 < p_4 < p_5$. Denote $\tau = \tau_{D \cap B(v, p_4)}$. For $x \in F$ we define
\begin{align*}
  f_1(x)
  & =
  \ex^x \expr{f(X(\tau)) \; ; \; X(\tau) \notin B(v, p_5)}
  \, , \\
  f_2(x)
  & =
  \ex^x \expr{f(X(\tau)) \; ; \; X(\tau) \in B(v, p_5)}
  \, .
\end{align*}
Clearly, $f = f_1 + f_2$, $f_1 = 0$ on $B(v, p_5) \setminus D$, $f_2 = 0$ on $(B(v, p_5))^c$, and both $f_1$ and $f_2$ are regular $\alpha$-harmonic in $D \cap B(v, p_4)$. We first estimate $f_1$.

By the strong Markov property, for $x \in D \cap B(v, p_3)$,
\begin{align*}
  f_1(x)
  & =
  \ex^x \expr{f(X(\tau)) \; ; \; X(\tau) \notin B(v, p_5)}
  \le
  \ex^x \expr{f(X(\tau_{B(v, p_4)})) \; ; \; X(\tau) \notin B(v, p_5)}
  \, .
\end{align*}
From~\eqref{eq:iw} it follows that
\begin{equation}
\label{eq:iw:estimate}
\begin{split}
  f_1(x)
  & \le
  c_5 \int_{B(v, p_5)^c} \int_{B(v, p_4)} \frac{G_{B(v, p_4)}(x, y) \, f(z)}{\varrho(y, z)^{d + \alpha}} \, \mu(dy) \, \mu(dz)
  \\ & \le
  \frac{c_5}{(1 - \frac{p_4}{p_5})^{d + \alpha}} \expr{\int_{B(v, p_4)} G_{B(v, p_4)}(x, y) \mu(dy)} \expr{\int_{B(v, p_5)^c} \frac{f(z)}{\varrho(v, z)^{d + \alpha}} \, \mu(dz)}
  \\ & \le
  \frac{c_5}{(1 - \frac{p_4}{p_5})^{d + \alpha}} \, \ex^x \tau_{B(v, p_5)} \, \Lambda_{v, p_3}(f)
  \, .
\end{split}
\end{equation}
Since $\ex^x \tau_{B(v, p_5)}$ is bounded, the upper bound for $f_1$ follows. It remains to estimate $f_2$.

Let $P = P_{v, p_4, p_3}$ be the function defined before the statement of the lemma. For $x \in D \cap B(v, p_3)$ we have
\begin{align*}
  f_2(x)
  & \le
  \int_{B(v, p_3)^c} P(x, y) \, f_2(y) \, \mu(dy)
  \le
  c_6 \, p_4^{-d} \int_{B(v, p_3)^c} f_2(y) \, \mu(dy)
  \, .
\end{align*}
Since $f_2(y) = 0$ for $y \in (B(v, p_5))^c$ and $f_2 \le f$, we conclude that
\begin{align*}
  f_2(x)
  & \le
  c_6 \, p_4^{-d} \int_{B(v, p_5) \setminus B(v, p_3)} f(y) \, \mu(dy)
  \le
  c_6 \, p_4^{-d} \Lambda_{v, p_3}(f)
  \, .
\end{align*}
This completes the proof.
\end{proof}

\begin{lemma}
\label{lem:factorization}
Suppose that $\alpha \in (0, 1)$ and $0 < p_1 < p_5 \le 1$. There are constants $c_8 = c_8(p_1, p_5, \alpha)$ and $c_8' = c_8'(p_1, p_5, \alpha)$ with the following property. If a nonnegative function $f$ is regular $\alpha$-harmonic in an open $D \sub B(v, p_5 2^{-m})$, where $v \in \cV_m$, and vanishes on $D^c \cap B(v, p_5 2^{-m})$, then
\begin{align}
\label{eq:factorization}
  c_8' \, \Lambda_{v, p_1 2^{-m}}(f) \, \ex^x \tau_D
  \leq
  f(x)
  & \leq
  c_8 \, \Lambda_{v, p_1 2^{-m}}(f) \, \ex^x \tau_D
  \, , &
  x \in D \cap B(v, p_1 2^{-m})
  \, .
\end{align}
\end{lemma} 

\begin{proof}
Again with no loss of generality we may assume that $m = 0$. Let $p_2, p_3$ satisfy $p_1 < p_2 < p_3 < p_5$, and let $\tau = \tau_{D \cap B(v, p_2)}$. We have $f = f_1 + f_2$, where
\begin{align*}
  f_1(x)
  & =
  \ex^x \expr{f(X(\tau)) \; ; \; X(\tau) \notin B(v, p_3)}
  \, , \\
  f_2(x)
  & =
  \ex^x \expr{f(X(\tau)) \; ; \; X(\tau) \in B(v, p_3)}
  \, .
\end{align*}
We estimate $f_1$ using~\eqref{eq:iw} as in~\eqref{eq:iw:estimate}, with $p_2$ and $p_3$ in place of $p_4$ and $p_5$. It follows that for $x \in D \cap B(v, p_1)$ we have
\begin{align*}
  f_1(x)
  & \le
  \frac{c_5}{(1 - \frac{p_2}{p_3})^{d + \alpha}} \, \ex_x \tau_{D \cap B(v, p_2)} \, \Lambda_{v, p_3}(f)
  \, .
\end{align*}
A similar lower bound holds with constant $c_5' (1 + \frac{p_2}{p_3})^{-d - \alpha}$. To estimate $f_2$, we use Lemmas~\ref{lem:escape} and~\ref{lem:upper}. For $x \in D \cap B(v, p_1)$,
\begin{align*}
  f_2(x)
  & \le
  \pr^x(X(\tau) \in B(v, p_3)) \, \sup_{y \in B(v, p_3)} f(y)
  \\ & \le
  \pr^x(X(\tau) \notin B(v, p_2)) \expr{c_7 \, \Lambda_{v, p_3}(f)}
  \\ & \le
  c_4 \, c_7 \, \Lambda_{v, p_3}(f) \, \ex^x \tau_{D \cap B(v, p_2)}
  \, .
\end{align*}
It follows that for some $C, C'$ dependent only on $p_j$ and $\alpha$,
\begin{align*}
  C' \, \Lambda_{v, p_3}(f) \, \ex^x \tau_{D \cap B(v, p_2)}
  \le
  f(x)
  & \le
  C \, \Lambda_{v, p_3}(f) \, \ex^x \tau_{D \cap B(v, p_2)}
  \, , &
  x & \in B(v, p_1)
  \, .
\end{align*}
Clearly $\ex^x \tau_{D \cap B(v, p_2)} \le \ex^x \tau_D$. The strong Markov property and Lemma~\ref{lem:escape} yield that for $x \in D \cap B(v, p_1)$ we also have
\begin{align*}
  \ex^x \tau_D
  & =
  \ex^x \tau_{D \cap B(v, p_2)} + \ex^x \expr{\ex^{X(\tau_{D \cap B(v, p_2)})}(\tau_D) \; ; \; X(\tau_{D \cap B(v, p_2)}) \notin B(v, p_2)}
  \\ & \le
  \ex^x \tau_{D \cap B(v, p_2)} + \pr^x \expr{X(\tau_{D \cap B(v, p_2)}) \notin B(v, p_2)} \sup_{y \in D} E^y \tau_D
  \\ & \le
  \expr{1 + c_4 \sup_{y \in B(v, 1)} E^y \tau_{B(v, 1)}} \ex^x \tau_{D \cap B(v, p_2)}
  \, .
\end{align*}
Obviously $\Lambda_{v, p_3}(f) \le \Lambda_{v, p_1}(f)$. On the other hand, by Lemma~\ref{lem:upper} for $x \in D \cap B(v, p_1)$ we have
\begin{align*}
  \Lambda_{v, p_1}(f)
  & \le
  \Lambda_{v, p_3}(f) + p_1^{-d - \alpha} \, \mu(D \cap B(v, p_3)) \, \sup_{y \in D \cap B(v, p_3)} f(x)
  \\ & \le
  \expr{1 + c_7 \, p_1^{-d - \alpha} \, \mu(B(v, 1))} \Lambda_{v, p_3}(f)
  \, .
\end{align*}
This proves~\eqref{eq:factorization}.
\end{proof}

\begin{proof}[Proof of Theorem~\ref{th:bhi}]
From Lemma~\ref{lem:factorization} with $p_1 = \frac{1}{2}$, $p_5 = \frac{\sqrt{3}}{2}$ we have for $x, y \in D \cap B(v, 2^{-m - 1})$
\begin{align*}
  f(x) \, g(y)
  & \le
  \expr{c_8 \, \Lambda_{v, \frac{1}{2}}(f) \, \ex^x \tau_{D'}} \expr{c_8 \, \Lambda_{v, \frac{1}{2}}(g) \, \ex^y \tau_{D'}}
  \\ & =
  \expr{\frac{c_8}{c_8'}}^2 \expr{c_8' \, \Lambda_{v, \frac{1}{2}}(f) \, \ex^y \tau_{D'}} \expr{c_8' \, \Lambda_{v, \frac{1}{2}}(g) \, \ex^x \tau_{D'}}
  \le
  \expr{\frac{c_8}{c_8'}}^2 \, f(y) \, g(x)
  \, ,
\end{align*}
where $D' = D \cap B(v, \frac{\sqrt{3}}{2})$.
\end{proof}

\noindent \textbf{Acknowlegements.}
The authors are deeply indebted to Andrzej St{\'o}s for discussion, advice and encouragement throughout the research and preparation of this article. Special thanks go to Krzysztof Bogdan for his careful reading of the manuscript and numerous suggestions.

\end{document}